\documentclass{article}
\usepackage{latexsym}
\usepackage{graphicx}
\usepackage{picfast}
\newtheorem{thm}{Theorem}
\newtheorem{conj}{Conjecture}
\newtheorem{prop}{Proposition}
\newtheorem{lem}{Lemma}
\newtheorem{defn}{Definition}
\newtheorem{quest}{Question}
\begin{document}
\begin{center}
{\bf THE INTERSECTION GRAPH CONJECTURE FOR LOOP DIAGRAMS}\\
\vspace{.2in}
{\footnotesize BLAKE MELLOR}\\
{\footnotesize Honors College}\\
{\footnotesize Florida Atlantic University}\\
{\footnotesize 5353 Parkside Drive}\\
{\footnotesize Jupiter, FL  33458}\\
{\footnotesize\it  bmellor@fau.edu}\\
\vspace{.2in}
{\footnotesize ABSTRACT}\\
{\ }\\
\parbox{4.5in}{\footnotesize \ \ \ \ \ The study of Vassiliev invariants for knots can be
reduced to the study of the algebra of chord diagrams modulo certain relations
(as done by Bar-Natan).  Chmutov, Duzhin and Lando defined the idea of the intersection
graph of a chord diagram, and conjectured that these graphs determine the
equivalence class of the chord diagrams.  They proved this conjecture in the case
when the intersection graph is a tree.  This paper extends their proof to the case
when the graph contains a single loop, and determines the size of the subalgebra
generated by the associated ``loop diagrams.''  While the conjecture is known to
be false in general, the extent to which it fails is still unclear, and this result
helps to answer that question.}\\
\vspace{1in}
\end{center}
\input{psfig.sty}
\section{Introduction} \label{S:intro}
In 1990, V.A. Vassiliev introduced the idea of {\it Vassiliev} or {\it finite type} knot
invariants, by looking at certain groups associated with the cohomology of the space of
knots.  Shortly thereafter, Birman and Lin~\cite{bn} gave a combinatorial description
of finite type invariants.  We will give a very brief overview of this combinatorial theory.
  For more details, see Bar-Natan~\cite{bn} and Chmutov, Duzhin and Lando~\cite{cdl1}.

We first note that we can extend any knot invariant to an invariant of
{\it singular} knots, where a singular knot is an immersion of $S^1$ in
3-space which is an embedding except for a finite number of isolated double
points.  Given a knot invariant $v$, we extend it via the relation:
$$\psfig{file=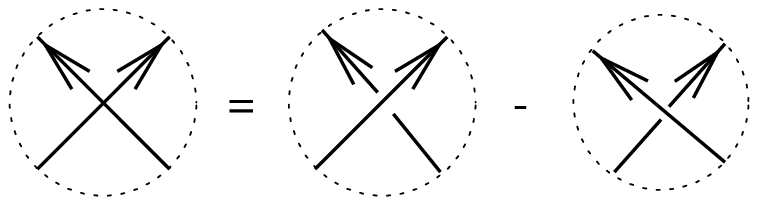}$$
An invariant $v$ of singular knots is then said to be of {\it finite type}, if
there is an integer $n$ such that $v$ is zero on any knot with more than $n$
double points.  $v$ is then said to be of {\it type} $n$.  The smallest such
$n$ is called the {\it order} of $v$.  We denote by $V_n$ the space generated
by finite type invariants of type $n$ (i.e., whose order is $\leq$ $n$).  We
can completely understand the space of finite type invariants by understanding
all of the vector spaces $V_n/V_{n-1}$.  An element of this vector space is
completely determined by its behavior on knots with exactly $n$ singular
points.  In addition, since such an element is zero on knots with more than
$n$ singular points, any other (non-singular) crossing of the knot can be
changed without affecting the value of the invariant.  This means that
elements of $V_n/V_{n-1}$ can be viewed as functionals on the space of
{\it chord diagrams}:
\begin{defn}
A {\bf chord diagram of degree n} is an oriented circle, together with $n$
chords of the circles, such that all of the $2n$ endpoints of the chords are
distinct.  The circle represents a knot, the endpoints of a chord represent
2 points identified by the immersion of this knot into 3-space.  The diagram
is determined by the order of the $2n$ endpoints.
\end{defn}
Note that since rotating a chord diagram does not change the order of the
endpoints of the chords, it does not change the diagram.  In particular,
turning a diagram ``upside down'' by rotating it 180 degrees gives the same
diagram.  We will use this fact later in the paper.

Functionals on the space of chord diagrams which are derived from knot
invariants will satisfy certain relations.  This leads us to the definition
of a {\it weight system}:
\begin{defn}
A {\bf weight system of degree n} is a linear functional $W$ on the space of
chord diagrams of degree $n$ (with values in an associative commutative ring
${\bf K}$ with unity) which satisfies 2 relations:
\begin{itemize}
    \item (1-term relation) $$\psfig{file=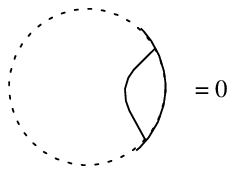}$$
    \item (4-term relation) $$\psfig{file=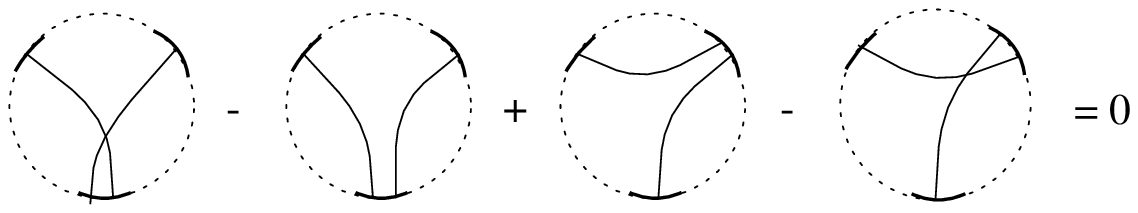}$$
    Outside of the solid arcs on the circle, the diagrams can be anything,
as long as it is the same for all four diagrams.
\end{itemize}
\end{defn}
It can be shown (see \cite{bl,bn,st}) that the space $W_n$ of weight
systems of degree $n$ is isomorphic to $V_n/V_{n-1}$.  For convenience,
we will usually take the dual approach, and simply study the space of
chord diagrams of degree $n$ modulo the 1-term and 4-term relations.
Bar-Natan~\cite{bn} and Kneissler~\cite{kn} have computed the dimension of these spaces for $n \leq
12$.  It is useful to combine all of these spaces into a graded
module via direct sum.  We can give this module a Hopf algebra
structure by defining an appropriate product and co-product:
\begin{itemize}
    \item  We define the product $D_1 \cdot D_2$ of two chord diagrams
$D_1$ and $D_2$ as their connect sum.  This is well-defined modulo the 4-term
relation (see \cite{bn}).
$$\psfig{file=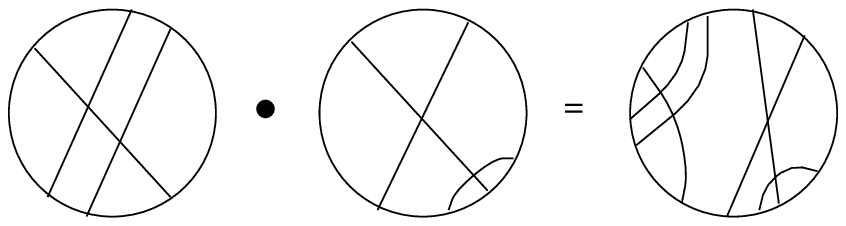}$$
    \item  We define the co-product $\Delta(D)$ of a chord diagram $D$ as
follows:
$$
\Delta(D) = {\sum_J D_J' \otimes D_J''}
$$
where $J$ is a subset of the set of chords of $D$, $D_J'$ is $D$ with all the
chords in $J$ removed, and $D_J''$ is $D$ with all the chords not in J
removed.
$$\psfig{file=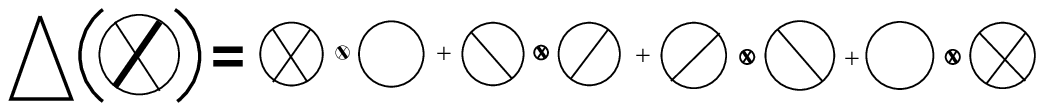}$$
\end{itemize}
It is easy to check the compatibility condition $\Delta(D_1\cdot D_2) = \Delta
(D_1)\cdot\Delta(D_2)$.  The rest of this paper is concerned with studying
parts of this Hopf algebra.

It is often useful to consider the Hopf algebra of bounded unitrivalent diagrams, rather
than chord diagrams.  These diagrams, introduced by Bar-Natan~\cite{bn} (there called
{\it Chinese Character Diagrams}), can
be thought of as a shorthand for writing certain linear combinations of chord
diagrams.  We define a {\it bounded unitrivalent diagram} to be a unitrivalent graph, with
oriented vertices, together with a bounding circle to which all the univalent vertices are
attached.  We also require that each component of the graph have at least one univalent
vertex (so every component is connected to the boundary circle).  We define the space $A$
of bounded unitrivalent diagrams as the quotient of the space of all bounded unitrivalent
graphs by the $STU$ relation:
$$\psfig{file=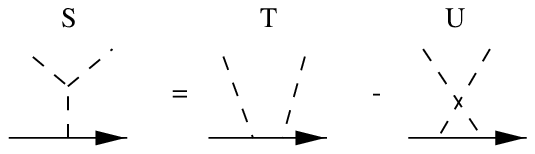}$$
As consequences of the $STU$ relation, the anti-symmetry ($AS$) and $IHX$ relations also hold in $A$:
$$\psfig{file=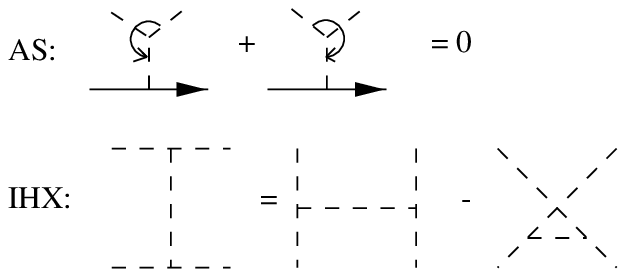}$$
Bar-Natan shows that $A$ is isomorphic to the algebra of chord diagrams.
\section{Intersection Graphs} \label{S:graphs}
\begin{defn}
Given a chord diagram $D$, we define its {\bf intersection graph} $\Gamma(D)$ as the graph
such that:
\begin{itemize}
    \item $\Gamma(D)$ has a vertex for each chord of $D$.
    \item Two vertices of $\Gamma(D)$ are connected by an edge if and only if the corresponding
chords in $D$ intersect, i.e. their endpoints on the bounding circle alternate.
\end{itemize}
\end{defn}
For example:
$$\psfig{file=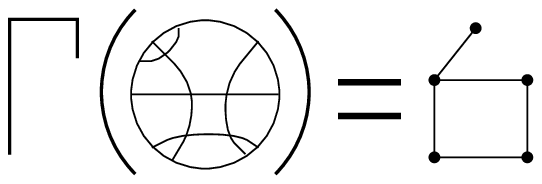}$$
Not all graphs can be the intersection graph for a chord diagram.  The simplest example of a graph
which cannot be an intersection graph occurs with 6 vertices:
$$\psfig{file=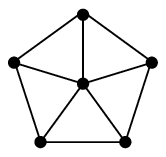}$$
Also, a graph can be the intersection graph for more than one chord diagram.  For example, there
are three different chord diagrams with the following intersection graph:
$$\psfig{file=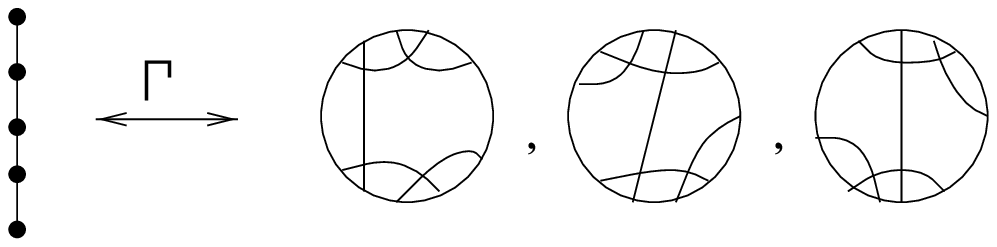}$$

However, these chord diagrams are all equivalent modulo the 4-term relation.  Chmutov,
Duzhin and Lando~\cite{cdl1} conjectured that intersection graphs actually determine
the chord diagram, up to the 4-term relation.  In other words, they proposed:
\begin{conj} \label{C:IGC}
If $D_1$ and $D_2$ are two chord diagrams with the same intersection graph, i.e. $\Gamma(D_1)$
= $\Gamma(D_2)$, then for any weight system $W$, $W(D_1) = W(D_2)$.
\end{conj}
This Intersection Graph Conjecture is now known to be false in general.  Morton
and Cromwell~\cite{mc} found a finite type invariant of type 11 which can distinguish some mutant knots,
and Le~\cite{le} and Chmutov and Duzhin~\cite{cd} have shown that mutant knots cannot be distinguished by
intersection graphs.  However, the conjecture is true in many special cases, and the exact extent to which
it fails is still unknown, and potentially very interesting.

The conjecture is known to hold in the following cases:
\begin{itemize}
    \item  For chord diagrams with 8 or fewer chords (checked via computer calculations);
    \item  For the weight systems coming from the defining representations of Lie
algebras {\it gl(N)} or {\it so(N)} as constructed by Bar-Natan in \cite{bn};
    \item  When $\Gamma(D_1) = \Gamma(D_2)$ is a tree (or, more generally, a linear
combination of forests) (see \cite{cdl1}).
\end{itemize}

The main result of this paper is to add one more case to the above list; namely, when the
intersection graph contains a single loop:
\begin{thm} \label{T:IGCloops}
If $D_1$ and $D_2$ are two chord diagrams such that $\Gamma(D_1)=\Gamma(D_2)=\Gamma$, and
$\Gamma$ has at most one loop, i.e. $\pi_1(\Gamma) = {\bf Z} {\rm \ or\ } 0$, then for any
weight system $W$, $W(D_1) = W(D_2)$.
\end{thm}
Our proof of this fact will closely follow the arguments of Chmutov et. al.~\cite{cdl2}, so we will begin by
recalling some definitions and results from that article, and then generalize them to prove our
result.
\section{Shares, Elementary Transformations, and Tree Diagrams} \label{S:shares}
We begin with the important idea of a {\it share} of a chord diagram:
\begin{defn}
Let $D$ be a chord diagram, and $C$ its collection of chords.  A {\bf share} is a subset
$S \subset C$ such that there exist four points $x_1,x_2,x_3,x_4$ in order along the circle so that:
\begin{itemize}
    \item  $x_i$ is not an endpoint of a chord;
    \item  For any chord $c \in S$, the endpoints of $c$ are in the arcs $x_1x_2$ and
$x_3x_4$;
    \item  For any chord $c \in C-S$, the endpoints of $c$ are in the arcs $x_2x_3$ and
$x_4x_1$.
\end{itemize}
In other words, we can divide the circle into 4 arcs so that no chord connects adjacent arcs.
\end{defn}
For example, in Figure~\ref{F:share} the chords contained in the dotted region form a share.  But
the three thick chords on their own do {\bf not} form a share, because there is a chord
which separates their endpoints.
    \begin{figure}
    \begin{center}
[share.pic]%    \input{share.pic}
    \caption{Share marked by dashed line} \label{F:share}
    \end{center}
    \end{figure}

Note that any single chord is a share, and the complement $C-S$ of a share $S$
is also a share.  There is one more important case:
\begin{lem} \label{L:bough}
Let $D$ be a chord diagram with $\Gamma(D)$ connected.  Suppose $D$ has a distinguished
chord $t$ (called the {\bf trunk}), and $v(t)$ is the vertex of $\Gamma(D)$
corresponding to $t$.  Then the chords of $D$ corresponding to the vertices
in a single component of $\Gamma(D)-v(t)$ form a share.  Such a share is
called a {\bf bough} of $t$.
\end{lem}
{\sc Proof:}  Let $S$ be a set of chords corresponding to a single component
of $\Gamma(D)-v(t)$.  Then no chord of $S^c = \Gamma(D)-\{S \cup v(t)\}$
intersects any chord of $S$.  In addition, since $S$ is connected, chords of
$S$ cannot be separated by chords of $S^c$; otherwise, a chord of $S^c$ would have
to intersect some chord on a path in $S$ connecting two chords of $S$.
  Now we can choose the minimal
arcs $x_1x_2$ and $x_3x_4$ containing all the endpoints of $S$.  There will be
two arcs, because $t$ intersects at least one chord of $S$, so the arcs must
be separated by the endpoints of $t$.  Therefore, there is a
chord of $S$ which connects the two arcs (namely, the one intersecting $t$).
$t$ will have one endpoint in $x_2x_3$ and one in $x_4x_1$, while all arcs
of $S^c$ will have both their endpoints in one of these two arcs.  Therefore,
$S$ is a share, and we are done.  $\Box$

Now we will consider a particular type of chord diagram called a {\it tree
diagram}:
\begin{defn}
A {\bf tree diagram} is a chord diagram $D$ whose intersection graph
$\Gamma(D)$ is a tree.
\end{defn}
The first question we ask is:  When do two tree diagrams have the same
intersection graph?  It is clear that we can permute the order of boughs along
some trunk of $D$ without changing the intersection graph.  For example:
$$\psfig{file=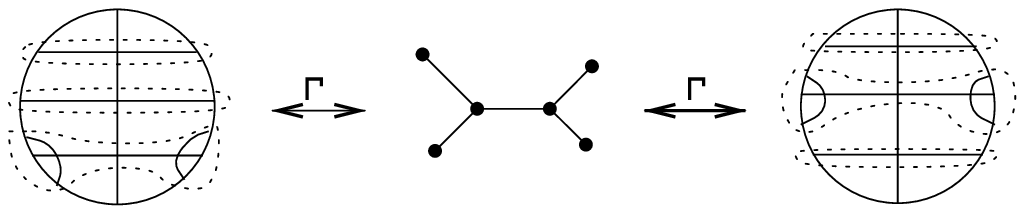}$$
In these diagrams the boughs correspond to the dotted regions.  We permute the
second and third boughs along the trunk.

We call such a permutation of boughs an {\it elementary transformation}.
Chmutov, Duzhin and Lando~\cite{cdl2} prove the following propositions:
\begin{prop} \label{P:tree_transform}
If $D_1$ and $D_2$ are tree diagrams such that $\Gamma(D_1) = \Gamma(D_2)$,
then $D_1$ can be transformed into $D_2$ by a sequence of elementary
transformations.
\end{prop}
\begin{prop} \label{P:tree_equivalent}
If $D_1$ and $D_2$ are tree diagrams which differ by an elementary
transformation, then $D_1$ and $D_2$ are equivalent modulo the 1-term and
4-term relations.
\end{prop}
Clearly, combining these two propositions gives us the Intersection Graph
Conjecture for tree diagrams.
The proof of the first proposition is straightforward (see \cite{cdl2}).  The
proof of the second proposition is much more complicated.  It
involves dividing the boughs along the trunk into ``upper boughs'' and
``lower boughs'', where the permutation affects solely the upper boughs, and
then inducting on the number of chords in the upper boughs.  Along the way,
Chmutov et al. prove several lemmas (see \cite{cdl2}):
\begin{lem} \label{L:Gen4term}
(The Generalized 4-term relation) $$\psfig{file=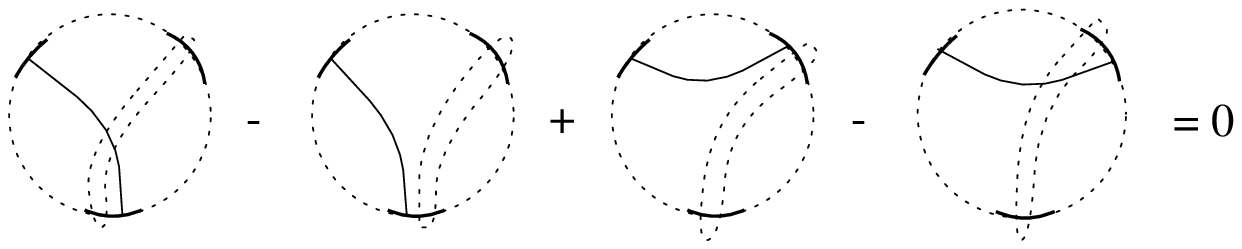}$$
As with the usual 4-term relation, the diagrams may be anything outside of the
solid arcs of the circle.
\end{lem}
Notice that the share is fixed while the chord moves around it.  The lemma
fails if we try to reverse these roles.  For example, if we apply the weight system
generated by the Kauffman polynomial (see Section~\ref{S:hopf}) to the linear combination of
chord diagrams in Figure~\ref{F:counterexample} we obtain $y^3x+4y^2x^2-yx^3 \neq 0$.
    \begin{figure}
    \begin{center}
[counterexample.pic]%    \input{counterexample.pic}
    \caption{Not a 4-term relation} \label{F:counterexample}
    \end{center}
    \end{figure}
\begin{lem} \label{L:corollary}
$$\psfig{file=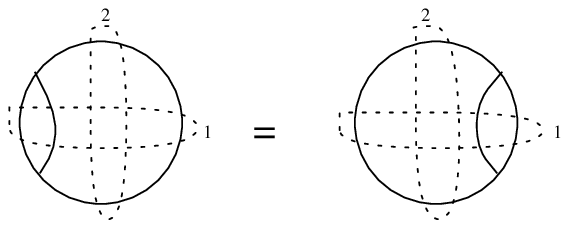}$$
\end{lem}
This lemma is a corollary of the generalized 4-term relation.  Notice that
neither of these lemmas make any assumptions about the diagram (i.e. it is
{\it not} necessarily a tree diagram).
\begin{lem} \label{L:decompose}
If D is the diagram on the left hand side below, and D-L is a tree, then:
$$\psfig{file=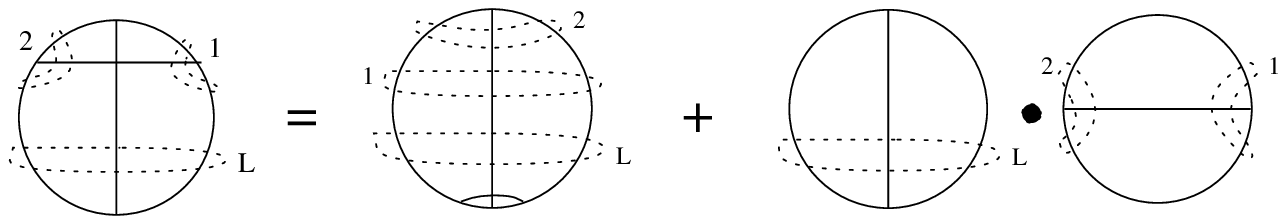}$$
\end{lem}
Chmutov et al. prove Lemma~\ref{L:decompose} only for tree diagrams, but a careful examination
of their proofs shows that we can weaken this hypothesis as above.  It
is only necessary that the {\it upper boughs} form a tree diagram.  The
components of $\Gamma(D)-v(t)$ corresponding to the lower boughs ($L$ in the
statements above) may contain loops.
\section{Loop Diagrams:  Definitions and Results} \label{S:loop}
\begin{figure}
    $$\psfig{file=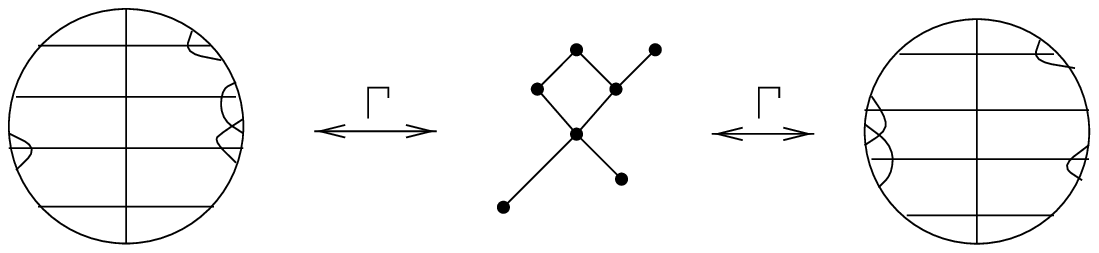}$$
    \caption{Reflecting a bough across a trunk} \label{F:reflect}
    \end{figure}

We closely model our treatment of loop diagrams on the treatment of tree
diagrams in \cite{cdl2}.
\begin{defn}
A chord diagram $D$ is called a {\bf loop diagram} if $\pi_1(\Gamma(D)) =
\bf Z$, i.e. the intersection graph has a single loop.
\end{defn}
Our first task is to determine the elementary transformations for loop
diagrams.  As with tree diagrams, we can certainly permute boughs along a
trunk.  However, these moves are not sufficient.  We can also reflect a
bough across a trunk, as shown in Figure~\ref{F:reflect}.
In the case of tree diagrams, this second move is a result of the first one,
but in the case of loop diagrams, there may be a special bough which
intersects the trunk twice, in which case the second move is independent of
the first one.  These two moves are now sufficient, so we define:
\begin{defn}
The {\bf elementary transformations} of a loop diagram are:
\begin{itemize}
    \item permuting boughs with respect to some trunk;
    \item reflecting a bough with respect to some trunk.
\end{itemize}
\end{defn}
Now we prove three propositions:
\begin{prop} \label{P:elementary}
If $D_1$ and $D_2$ are loop diagrams with $\Gamma(D_1) = \Gamma(D_2)$, then
$D_2$ can be obtained from $D_1$ by a sequence of elementary transformations.
\end{prop}
\begin{prop} \label{P:permutation}
If loop diagrams $D_1$ and $D_2$ differ by a permutation of boughs, then
$D_1$ and $D_2$ are equivalent modulo 1-term and 4-term relations.
\end{prop}
\begin{prop} \label{P:reflection}
If loop diagrams $D_1$ and $D_2$ differ by a reflection of a bough, then
$D_1$ and $D_2$ are equivalent modulo 1-term and 4-term relations.
\end{prop}
Clearly, combining these three propositions and the Intersection Graph Conjecture
for tree diagrams will give us Theorem 1.
\section{Proof of Intersection Graph Conjecture for loop diagrams} \label{S:proof}
\begin{figure}
    $$\psfig{file=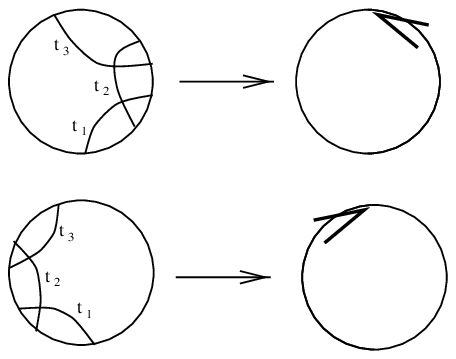}$$
    \caption{Order of chords induces orientation} \label{F:orient}
    \end{figure}

{\sc Proof of Proposition~\ref{P:elementary}:}  Let $\Gamma$ denote the intersection graph
$\Gamma(D_1) = \Gamma(D_2)$.  $\Gamma$ has a unique minimal loop of length
$\geq$ 3.  Choose a vertex $v_1$ on the loop, and a direction around the loop,
and use these choices to order the other vertices on the loop $v_1, v_2,
\dots, v_n$.  Now let $t_i$ and $t_i'$ be the chords corresponding to $v_i$
in $D_1$ and $D_2$ respectively.  The ordering of the vertices gives an
ordering of these chords in each diagram.  If $n > 3$, this ordering induces an orientation on
the bounding circles of the two diagrams, as shown in Figure~\ref{F:orient}.
If $n = 3$, we also induce an orientation on the bounding circle, though it's not so easy to see.
In this case we have three chords $t_i, t_2, t_3$, and the order of their endpoints moving clockwise
around the bounding circle is either $123123$ or $132132$ (see Figure~\ref{F:n=3}).
    \begin{figure}
    \begin{center}
[n3.pic]%    \input{n3.pic}
    \caption{The case when n=3} \label{F:n=3}
    \end{center}
    \end{figure}
In the first case, we will say the chords induce a clockwise orientation on the bounding circle,
and in the second case the chords induce a counterclockwise orientation.

\begin{figure}
    $$\psfig{file=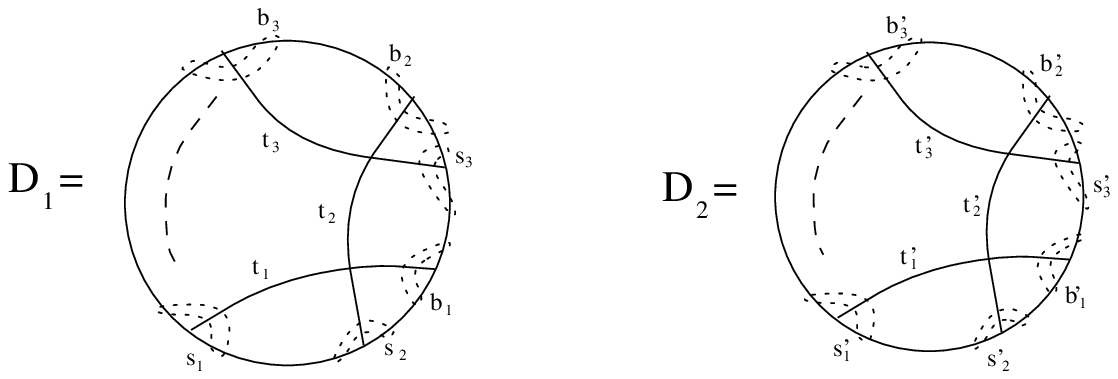}$$
    \caption{Counterclockwise orientation} \label{F:counter}
    \end{figure}

If these induced orientations do not agree, we can reflect the bough of $t_1$
in $D_1$ containing $\{t_2,\dots,t_n\}$ across $t_1$, which will reverse the
induced orientation.  So, via this elementary transformations, we may assume
that $D_1$ and $D_2$ are both oriented counterclockwise, as shown in Figure~\ref{F:counter}.

Now we can permute the boughs of $t_i$ so their order corresponds to the order
of the boughs along $t_i'$.  Since the other boughs of $t_i$ only intersect
$t_i$ once (the intersection graph has only one loop), each has a
distinguished trunk.  We can then permute boughs along these trunks.  As we
continue inductively, all further boughs will have a distinguished trunk, so
we can complete the transformation of $D_1$ to $D_2$ via permutations of
boughs.  $\Box$

\begin{figure}
    $$\psfig{file=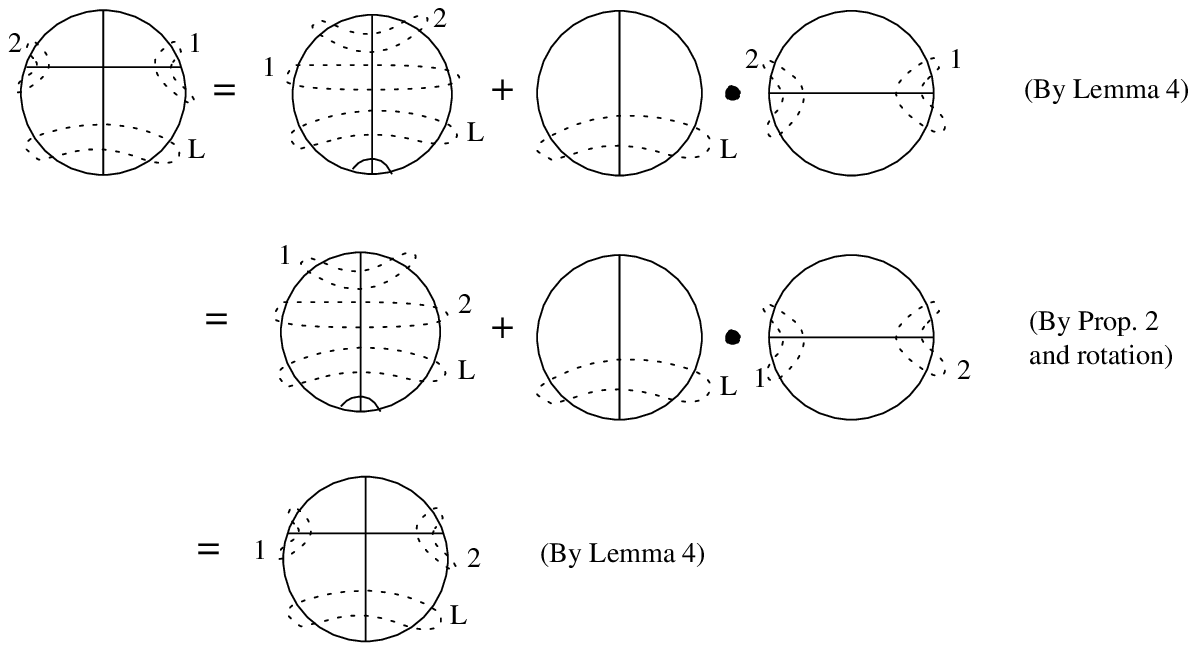}$$
    \caption{Proof of Lemma~\ref{L:rotate}} \label{F:proof}
    \end{figure}

To prove Proposition~\ref{P:permutation}, we will need the following lemma, which allows us to
rotate boughs in our chord diagrams:
\begin{lem} \label{L:rotate}
$$\psfig{file=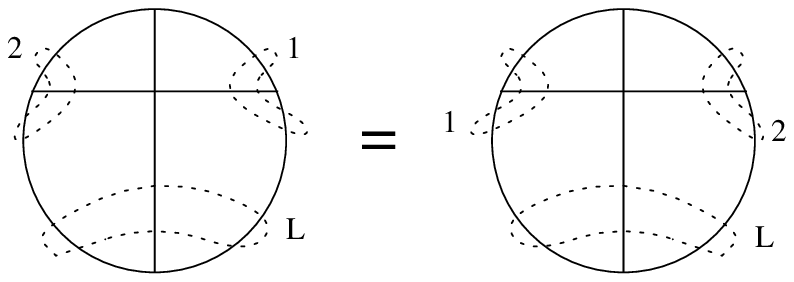}$$
where $L$ is any share, but $D-L$ is a tree diagram.  Here the shares 1 and 2 have each been
rotated 180 degrees (not reflected).
\end{lem}
{\sc Proof:}  By Lemma~\ref{L:decompose}, keeping in mind our observation that $L$ can be any
share, we get the equalities in Figure~\ref{F:proof}, which prove the lemma. $\Box$

{\sc Proof of Propostion~\ref{P:permutation}:}  First, we consider the case when one of the
boughs being permuted contains the loop; i.e. we show:
$$\psfig{file=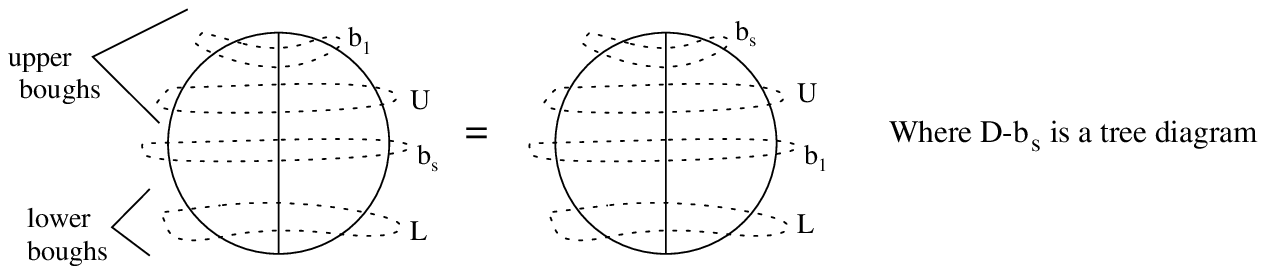}$$
Following Chmutov et. al.~\cite{cdl2} we show this by induction on the {\it complexity} of the
diagrams.  If we permute the boughs of $D$ by a permutation $\pi$, the
{\it lower boughs} are the boughs below $b_s$ in the diagram, and the
{\it upper boughs} are the boughs above $b_s$.  Then the
{\it complexity} $c(D,\pi)$ is the total number of chords in the upper boughs.

{\sc Base Case:}  When $c(D,\pi)$ = 1, then $b_1$ is just one chord.  So we
can move $b_1$ past $b_s$ via Lemma~\ref{L:corollary}.  Here share 1 of Lemma~\ref{L:corollary} is
$L \cup trunk$, and share 2 is $b_s$.

{\sc Inductive Case:}  Assume the proposition is true for $c(D,\pi) < m$.  We
will show it holds for $c(D,\pi) = m$.  The proof is by a chain of equalities
of chord diagrams:\\
$\psfig{file=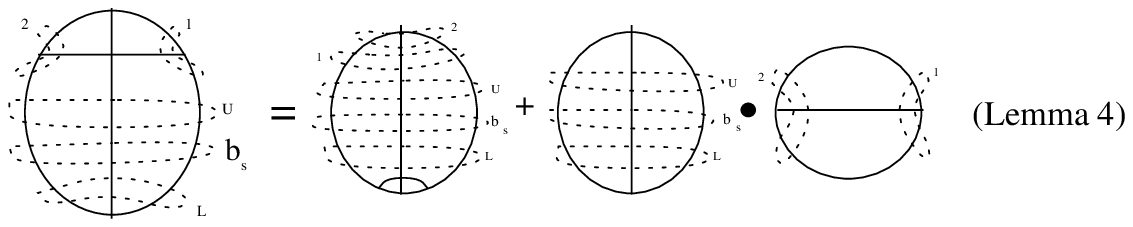}$\\
$\psfig{file=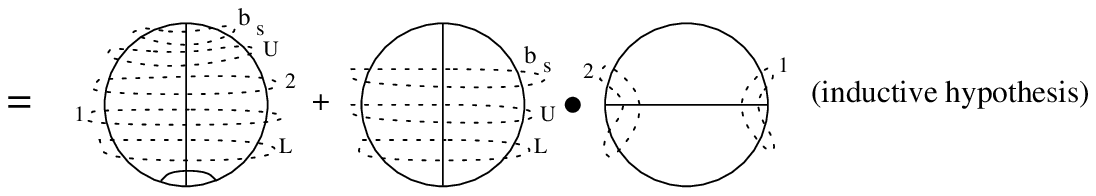}$\\
$\psfig{file=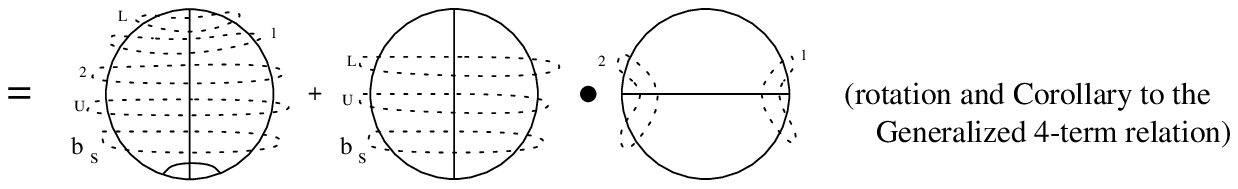}$\\
$\psfig{file=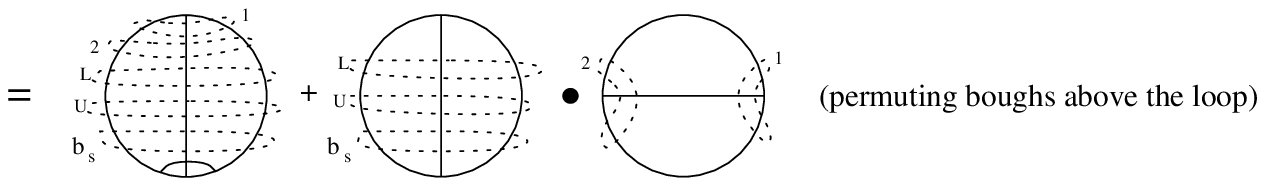}$\\
$\psfig{file=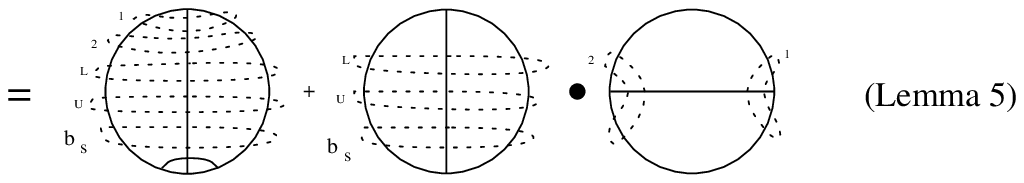}$\\
$\psfig{file=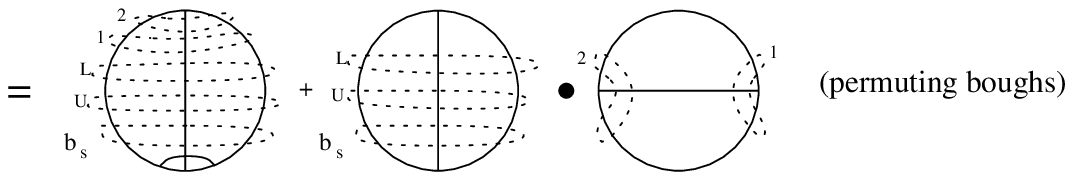}$\\
$\psfig{file=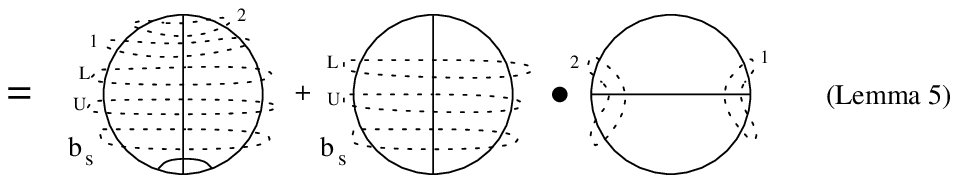}$\\
$\psfig{file=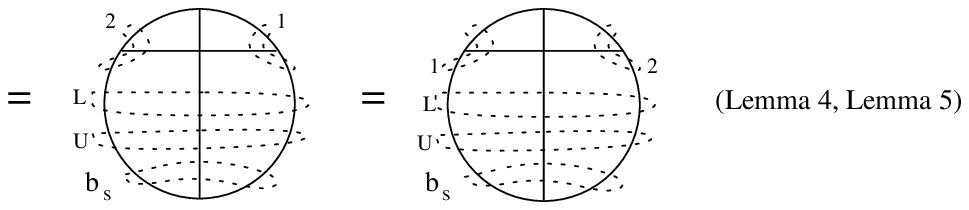}$\\
$\psfig{file=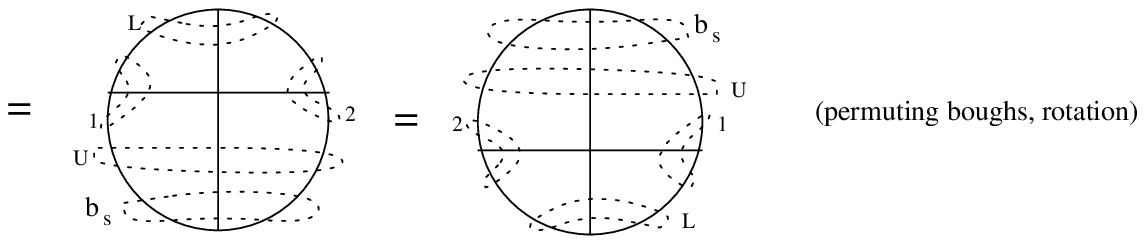}$\\
Finally, we consider the case when neither of the boughs being permuted
contains the loop.
If $L$ contains the loop, then we are done by Proposition~\ref{P:tree_equivalent} (using the
weakened hypothesis).  If $U$ contains the loop, we simply rotate the diagram
180 degrees, permute the boughs using Proposition~\ref{P:tree_equivalent}, and then rotate back.

This completes the proof of Propositon~\ref{P:permutation}. $\Box$

\begin{figure}
    $$\psfig{file=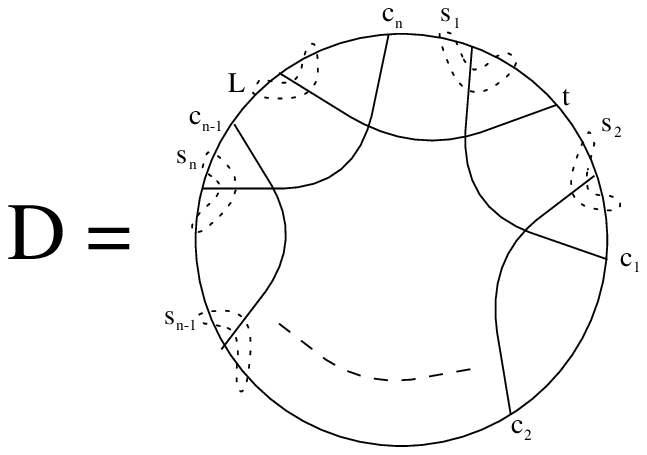}$$
    \caption{Normal Form} \label{F:normal}
    \end{figure}

Before we prove our final proposition, we will prove one more convenient lemma:
\begin{lem} \label{L:transfer}
$$\psfig{file=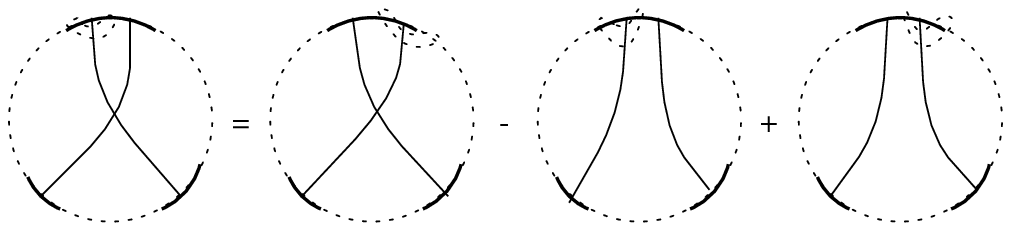}$$
\end{lem}
{\sc Proof:}  The proof is simply an application of the Generalized 4-term
relation:
$$\psfig{file=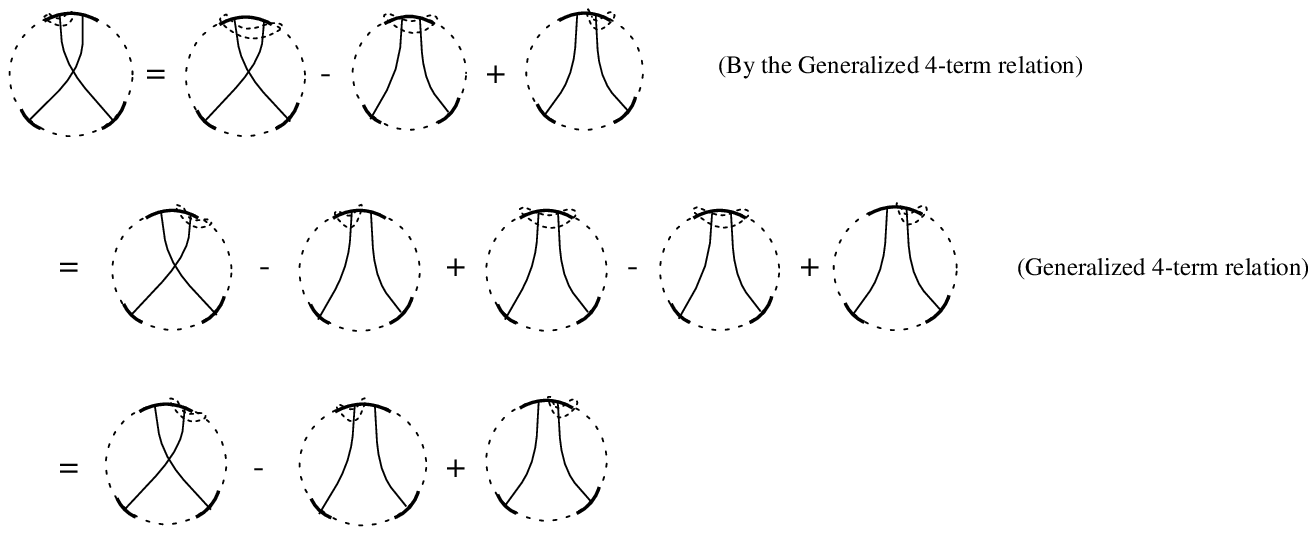}\Box$$
Now we can complete our argument by proving the final proposition:

{\sc Proof of Proposition~\ref{P:reflection}:}  We first consider the case when the bough does
not contain the loop (and so intersects the trunk only once).  In this case,
reflecting the bough is just the result of permuting its sub-boughs to
reverse their order (and doing the same on lower levels as necessary).  So
this case is a consequence of Proposition~\ref{P:permutation}.  The next case is when the bough
does contain the loop, but only intersects the trunk once.  In this case, we
can either accomplish the reflection by permuting boughs, or we reach a
stage when we wish to reflect the loop across a chord which it intersects
twice.  So we are reduced to the case of reflecting the loop across a trunk
t which it intersects twice.  By permuting boughs, we can put the diagram into a "normal form",
as shown in Figure~\ref{F:normal}.
Again, the proof will work by induction.  In this case, we will induct on
$max\{n | |s_n| > 0\}$, where the $s_n$'s are numbered clockwise as shown,
with $L = s_0$, and $|s_n|$ is the number of chords in the share $s_n$.  We
will call this the {\it heft} of the diagram.

\begin{figure}
    $$\psfig{file=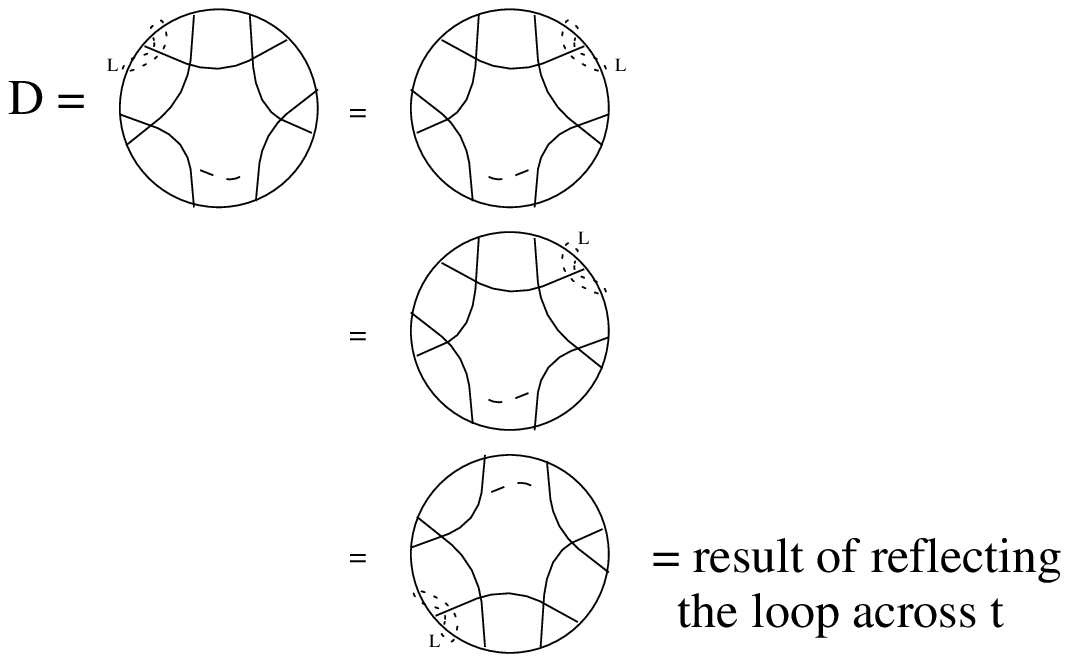}$$
    \caption{Proof of the Base Case} \label{F:base_pf}
    \end{figure}

{\sc Base Case:}  heft = 0.  This case is easily proved using Proposition~\ref{P:permutation}
and Lemma~\ref{L:transfer}.  See Figure~\ref{F:base_pf}.

{\sc Induction:}  Assume Proposition~\ref{P:reflection} holds when the heft is less than
{\it m}, and that the diagram $D$ has heft $m$.  Then Lemma~\ref{L:transfer} tells us that:
$$\psfig{file=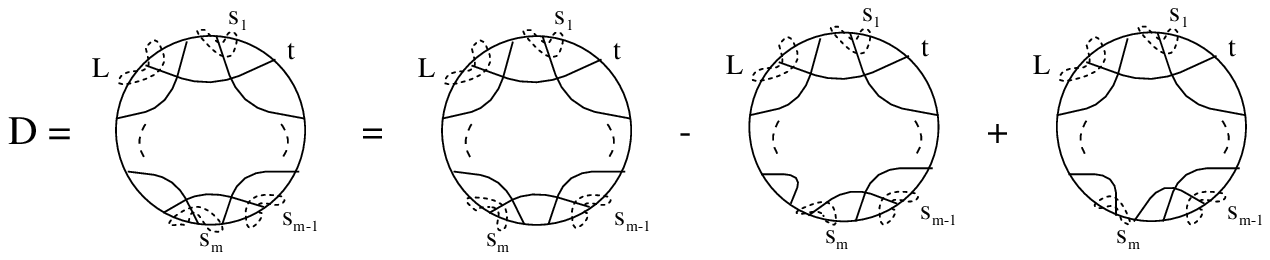}$$
The last two diagrams on the left hand side are trees, so we can reflect the
boughs through the trunk via permutations.  The first diagram on the left hand
side can be rewritten, permuting boughs, as:
$$\psfig{file=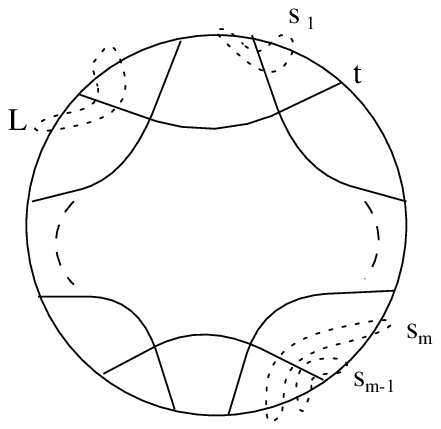}$$
Therefore, this diagram has heft $m-1$, and so the bough can be reflected by
our inductive hypothesis.  Doing the reflection on these three diagrams, and
then recombining them by Lemma 6, gives us the reflection in D:
$$\psfig{file=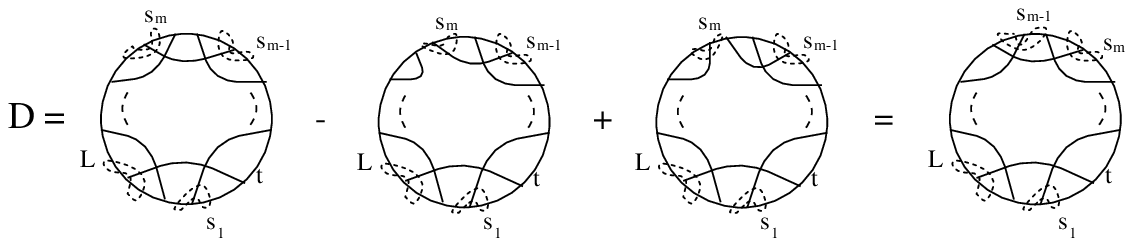}$$
This completes the proof of Proposition~\ref{P:reflection}, and hence of Theorem~\ref{T:IGCloops}.  $\Box$
\section{Hopf Algebra of Tree and Loop Diagrams} \label{S:hopf}
Chmutov, Duzhin and Lando~\cite{cdl3} use the fact that tree diagrams are
determined by their intersection graphs to compute the subalgebra of the
Hopf algebra of chord diagrams which is generated by the tree diagrams.  They
denote this subalgebra (the {\it forest subalgebra}) by {\it A}.  To be
precise, they prove:
\begin{thm} \label{T:tree_algebra}
The Hopf algebra {\it A} is isomorphic to the polynomial
algebra ${\bf Q}[x_1,x_2,...]$, where the grading of every $x_n$ is $n$.
\end{thm}
In this section we will make the analogous computation for the algebra
generated by both tree and loop diagrams.  We will rather unimaginatively
call this algebra the {\it loop algebra}, and denote it by {\it B}.  Our goal
in this section is to prove:
\begin{thm} \label{T:loop_algebra}
The Hopf algebra {\it B} is isomorphic to the polynomial algebra
$${\bf Q}[x_1,x_2,x_3,x_4^1,x_4^2,x_5^1,...,x_n^1,...,x_n^{n-2},...]$$
where the grading of every $x_n^i$ is $n$, and the primitive space of grading
$n$ has dimension $1$ if $n \leq 3$, and $n - 2$ if $n \geq 3$.
\end{thm}

The Hopf algebra will certainly be isomorphic to some such polynomial
algebra.  The content of the theorem is in computing the exact dimension of
the primitive space in each grading.  We do this by first finding an upper
bound for this dimension, and then finding sufficiently many independent
primitive elements to show that this upper bound is in fact the dimension of
the primitive space.  The first half of this program is accomplished via the
following proposition:
\begin{prop} \label{P:graph4term}
The 4-term relation on {\it B} induces three relations on intersection
graphs with at most one loop, shown in Figure~\ref{F:graph4term}.
    \begin{figure}
    $$\psfig{file=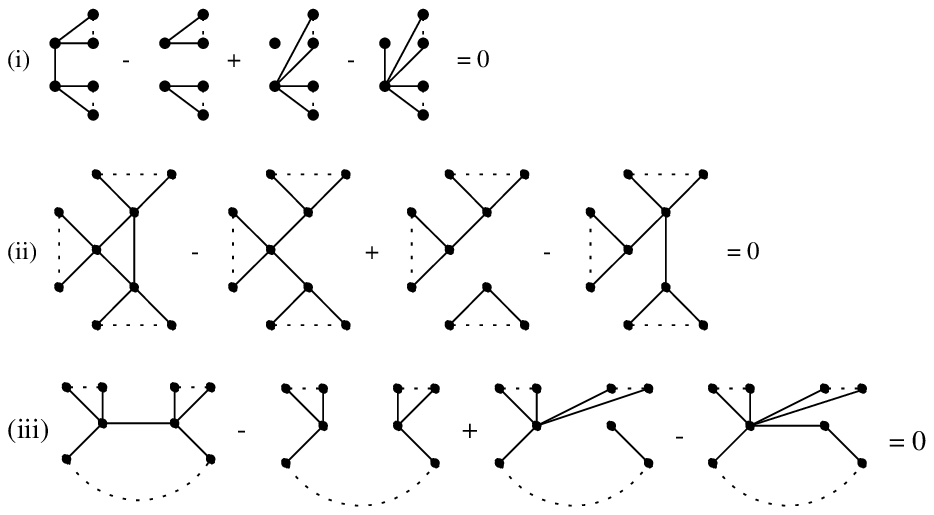}$$
    \caption{4-term relation for intersection graphs} \label{F:graph4term}
    \end{figure}
In these pictures, the graphs are identical outside of the region shown, and the dotted lines
indicate a group of edges incident to a single vertex.  In the third relation, the dotted
line along the bottom indicates that there is another path in the graph connecting the two vertices.
\end{prop}
{\sc Proof:}  The key to this proposition is Theorem~\ref{T:IGCloops}.  This allows us to say
that if the equalities hold for {\it any} chord diagrams with the given
intersection graphs, then they will hold for {\it all} such diagrams, since
the diagrams will be equivalent by Theorem~\ref{T:IGCloops}.  Then we will have the desired
relations induced on intersection graphs.  So we prove the equalities by
choosing nice chord diagrams for which the proofs are easy.

The first relation is induced by Lemma~\ref{L:decompose} (with our weakened hypotheses), with
an extra term which arises if we repeat the proof without using the 1-term
relation (see \cite{cdl2}).

The second relation results from the following equalities of chord diagrams:
$$\psfig{file=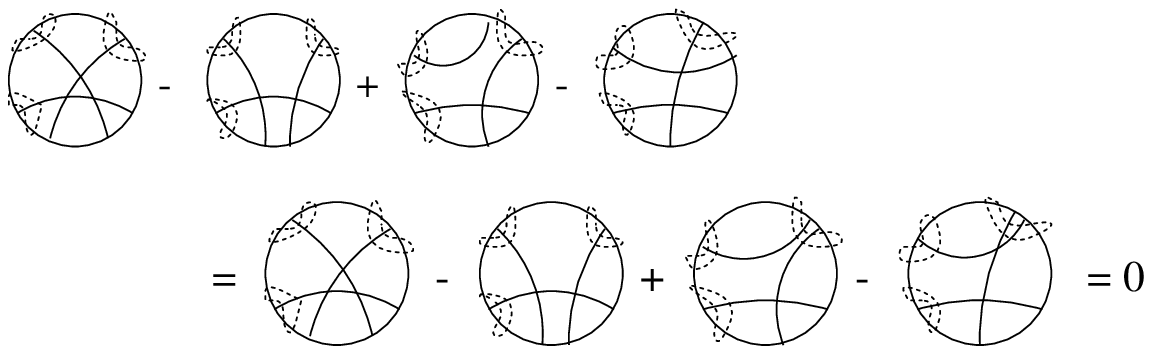}$$

The last relation is the most complicated to prove.  We first consider the
following two 4-term relations:
$$\psfig{file=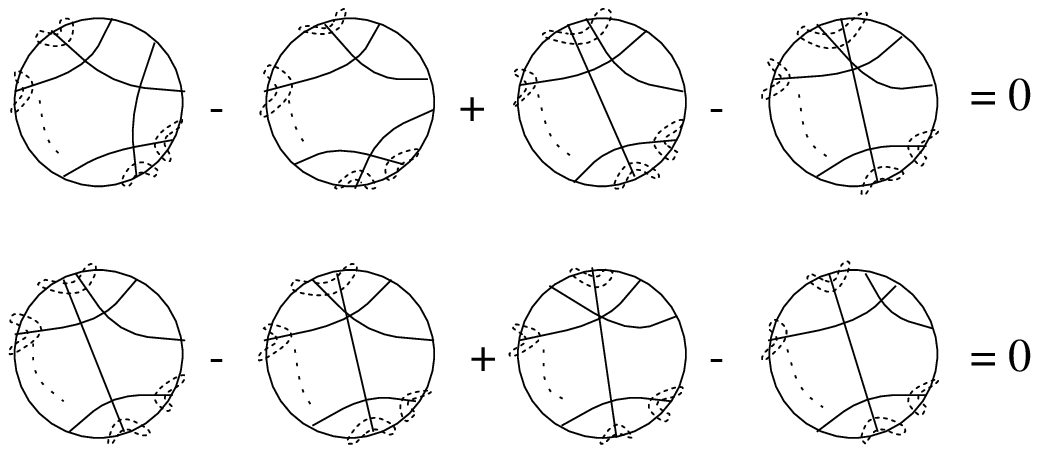}$$

By Theorem~\ref{T:IGCloops}, the last two terms on the left-hand side of the second equation
can be rewritten as:
$$\psfig{file=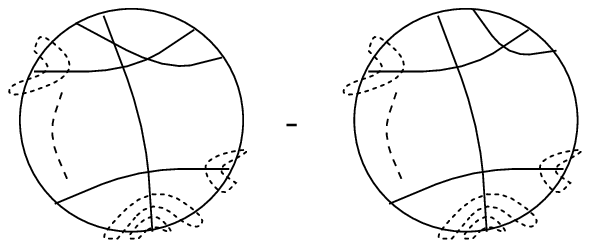}$$
So we consider a third 4-term relation:
$$\psfig{file=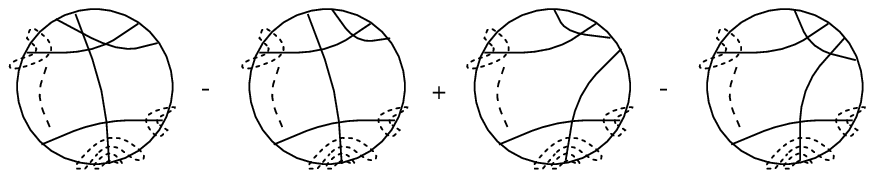}$$
Subtracting the second 4-term relation from the sum of the other two, we get:
$$\psfig{file=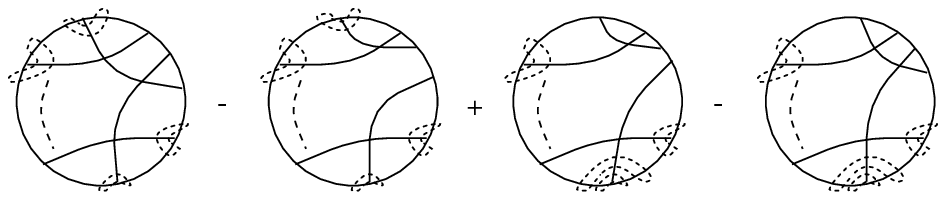}$$
which induces the desired relation on intersection graphs.  $\Box$

This proposition tells us that any tree with {\it n} vertices is equivalent,
modulo decomposable elements, to $\psfig{file=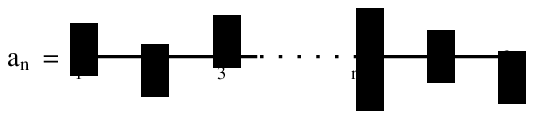}$.  Similarly, any graph
with {\it n} vertices whose only loop is a triangle is equivalent to $2a_n$,
and any loop graph with {\it n} vertices and a loop of length {\it k} is
equivalent to $\psfig{file=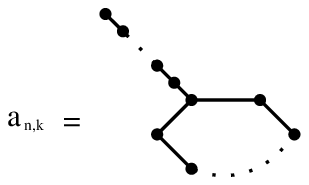}$.  This last equivalence is because the
middle two terms of part (iii) above cancel modulo decomposable elements
(they are both trees with {\it n} vertices), and then by repeated use of part
(i).  Inductively, we see that $\{a_n, a_{n,k}\}$ are generators for the
algebra of chord diagrams (though we don't yet know they are independent).
We would like to find a set of independent {\it primitive} generators in
each degree.  Recall that an element {\it p} of a Hopf algebra is
{\it primitive} if $\Delta(p) = 1\otimes p + p\otimes 1$, so no decomposable
element can be primitive.  The set of independent primitive elements can be no larger
than the set of graphs $\{a_n, a_{n,k}\}$ described above, so we have placed
an upper limit on the dimension of the primitive space of grading {\it n} of
the intersection graphs of loop diagrams, and hence (by Theorem~\ref{T:IGCloops}) on the
dimension of the primitive space of grading {\it n} of {\it B}.  This upper limit is 1 when
$n\leq 3$, and $n-2$ when $n\geq 3$.  So our goal is to show that these
upper limits are in fact the dimensions of the spaces by exhibiting a
sufficient number of primitive elements.
\begin{defn}
$p_n = \sum_J{(-1)^{|J|}a_{n,J}}$, where $J$ is some subset of the edges of
$a_n$, and $a_{n,J} = a_n-J$.
$p_{n,k} = \sum_J{(-1)^{|J|}a_{n,k,J}}$.
\end{defn}

We can show that all the elements $\{p_n, p_{n,k}\}$ are primitive directly, but
it is more elegant to use bounded unitrivalent diagrams.  Bar-Natan~\cite{bn} has shown
that an element of {\it A} is primitive if and only if it is a linear
combination of bounded unitrivalent diagrams with {\it connected} interiors (i.e. the
diagram minus its boundary circle is connected).  So it suffices to show that
$\{p_n, p_{n,k}\}$ have this form.  It is useful to recall the ``wheel with {\it k}
spokes,'' introduced by Chmutov and Varchenko~\cite{cv}:
$$\psfig{file=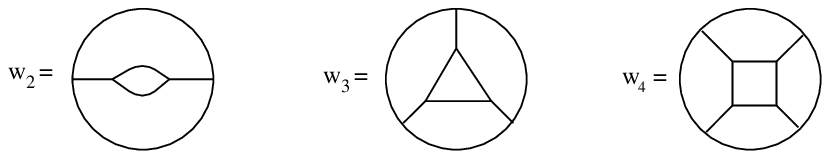}$$
\begin{prop} \label{P:primitive}
$$\psfig{file=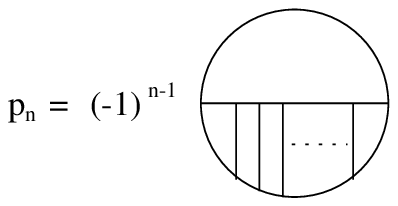}$$
$$\psfig{file=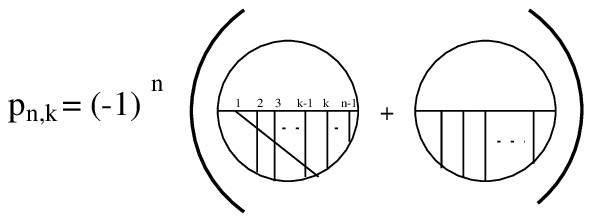}$$
$$p_{n,n} = (-1)^nw_n$$
Therefore, all the elements $p_n$ and $p_{n,k}$ are primitive.
\end{prop}
{\sc Proof:}  The first part of the proposition was noticed by Chmutov and Varchenko~\cite{cv}.
To prove it, we pick an edge {\it e} of $a_n$ and rewrite:
$$p_n = \sum_{J\ni e}{(-1)^{|J|}(a_{n,J}-a_{n,J-e})}$$
Each term of this sum is an STU relation, so we are left with a sum of bounded unitrivalent diagrams which
have inserted a ``leg'' in place of the edge {\it e} of $a_n$; there are half as many
terms in this sum as in the original one.  We continue this process with each edge in
turn, eventually obtaining $p_n = (-1)^{|all\ edges|}D = (-1)^{n-1}D$, where $D$ is the
bounded unitrivalent diagram in the proposition.  Via exactly the same argument, we find that $p_{n,n} =
(-1)^nw_n$.  The proof for $p_{n,k}$ is the same, except that we do not consider the
edge of $a_{n,k}$ where the loop reattaches, so we end up with 2 terms instead of one.
In this case, the final factor of -1 (to give a coefficient of $(-1)^n$ rather
than $(-1)^{n-1}$) comes from the anti-symmetry relation.  $\Box$

It only remains to show that $\{p_n,p_{n,k}\}_{k=4}^n$ is independent.
Chmutov, Duzhin and Lando~\cite{cdl3} show that the $p_n$'s are
non-zero using weighted graphs, but this
approach seems difficult to generalize.  Instead, we take a more direct
approach, and find a weight system $W$ such that $\{W(p_n),W(p_{n,k})\}$ is
independent.  We will use the weight system which arises from the Kauffman
polynomial, as described by Meng~\cite{me}.  We say that a knot is
{\it semi-oriented} if it is continuously oriented except at a finite number
of points.  Next, we recall the defining
skein relation for the Kauffman polynomial~\cite{ka}.  We modify the
relation slightly to give an invariant, up to {\it ambient} isotopy, of
semi-oriented knots.  We define links $L+, L-, L*, L\# and L!$ as in
Figure~\ref{F:skein} (where the links are the same outside of the region shown).
    \begin{figure}
\begin{center}
[skein-links.pic]%    $$\input{skein_links.pic}$$
    \caption{Links of the skein relation} \label{F:skein}
\end{center}
    \end{figure}
\begin{defn}
The Kauffman polynomial $F$ is the invariant of links, up to ambient isotopy,
defined by the skein relation:
$$bF(L+) - b^{-1}F(L-) = v(F(L\#) - F(L*))$$
In particular, this means that:
$$F(0 |) = \left({1+\frac{b^{-1}-b}{v}}\right)F(|)$$
\end{defn}

To obtain the weight system used by Meng, we make the following substitutions:
$$b = exp(-{1\over 2}yh);\ v = exp(-{1\over 2}xh) - exp({1\over 2}xh)$$
This immediately gives us the formula:
$$F(L!) = (exp({1\over 2}(y-x)h) - exp({1\over 2}(y+x)h))(F(L\#) - F(L*)) + (exp(yh) - 1)F(L-)$$

It is clear that the coefficient of $h^n$ is a finite type invariant of type
$n$.  To compute a weight system $W$ associated to the invariant, we isolate the
first non-zero coefficient, obtaining the relations:
$$W(L!) = yW(L-) + xW(L*) - xW(L\#)$$
$$W(O|) = (1-\frac{y}{x})W(|)$$
 
Now we want to evaluate this weight system on our primitive elements.  As an
example, we will compute $W(p_2)$ explicitly.  Note that dots on the boundary circles
indicate where the orientation reverses (so they always arise in pairs).  If two such
dots are connected by an arc containing no endpoints of chords, they can be moved
together by an isotopy of the diagram and cancelled.
$$W(p_2) = W\left(\psfig{file=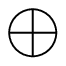} - \psfig{file=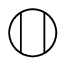}\right) =
W\left(\psfig{file=cross.ps}\right)$$
$$W\left(\psfig{file=cross.ps}\right) = yW\left(\psfig{file=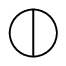}\right) +
xW\left(\psfig{file=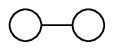}\right) - xW\left(\psfig{file=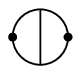}\right)$$
$$= x(yW(OO)+xW(O)-xW(O)) - x(yW(O)+xW(O)-xW(OO))$$
$$= x(y(1-\frac{y}{x})) - x(y+x-x(1-\frac{y}{x}))$$
$$= y(x-y) - x(2y) = -y(x+y)$$
\begin{lem} \label{L:weights}
$W(p_n) = -y(x+y)^{n-1}$ for all n.
$W(p_{n,k}) = (x+y)^{n-k}W(p_{k,k})$ for all $k<n$.
\end{lem}
{\sc Proof:}  We prove this lemma by induction.  We begin with the first statement.
First, we recall that (by Theorem 1):
$$\psfig{file=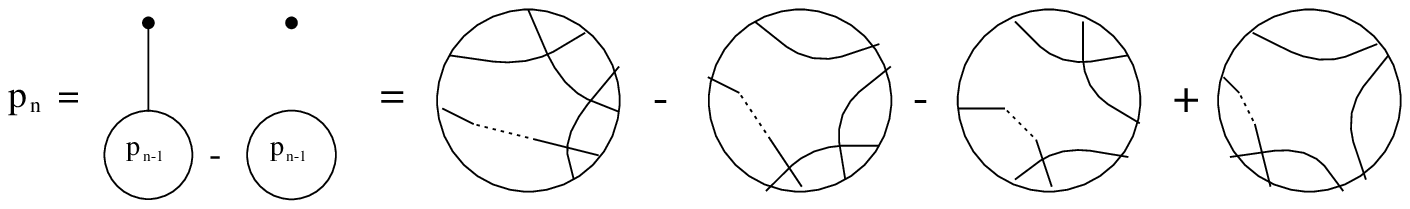}$$
The second and fourth terms disappear by the 1-term relation.  From the remaining terms,
we see (for $n>2$):
$$W(p_n) = yW(p_{n-1}) + xW\left(\psfig{file=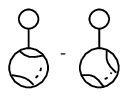}\right) - xW(p^*_{n-1}) =
yW(p_{n-1}) - xW(p^*_{n-1})$$
Where $p^*_n$ is defined by:
$$\psfig{file=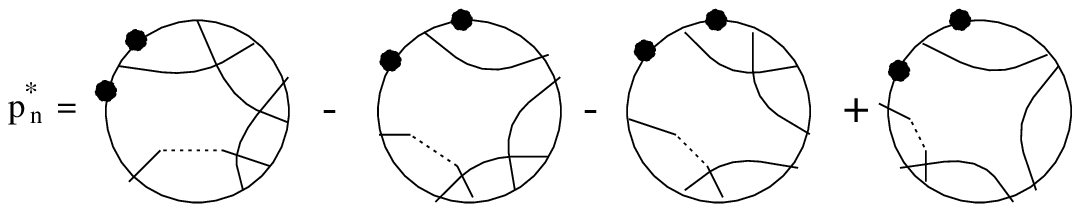}$$
The remainder of the diagrams are the same as $p_n$.

Similarly, noting that $W\left(\psfig{file=theta+.ps}\right) = 2y$, and that the changes
of orientation in $p^*_n$ effectively reverse the signs of the last two terms of the
weight system relation, we find that (for $n>2$):
$$W(p^*_n) = -yW(p_{n-1}) + xW(p^*_{n-1}) = -W(p_n)$$
A direct computation shows that $W(p^*_2) = -W(p_2)$.  Therefore,
$$W(p_n) = (x+y)W(p_{n-1})$$
A simple induction, together with our computation of $W(p_2)$, then gives the result.

The second statement is proved similarly, except that for the base case $p_{k+1,k}$, some
additional analysis of $W(p_{k,k})$ is required.  $\Box$
\begin{lem} \label{L:circle}
$W(p_{n,n}) = -2W(p_n)+yW(p_{n-1})+2x^2W(p_{n-2})+3x^2W(p_{n-2,n-2})-2x^3W(p_{n-3,n-3})$
for $n>6$.
\end{lem}
{\sc Proof:}  The proof of this lemma, while long, is elementary.  It is similar in concept
to the proof of Lemma~\ref{L:weights}, and is left to the industrious reader.  $\Box$

Meng~\cite{me} noted that for any chord diagram $D$, $W(D)$ has a factor $y(x+y)$.
Together with our results above, this implies that any $W(p_n),\ n>2$ or
$W(p_{n,k}),\ n>k$ has a factor $(x+y)^2$.  However, this is not the case for
$W(p_{n,n})$:
\begin{lem} \label{L:factor}
$W(p_{n,n})$ is not divisible by $(x+y)^2$.
\end{lem}
{\sc Proof:}  For n = 4,5,6 we show this by direct computation.
$$W(p_{4,4}) = y(x+y)(6x^2 + 3xy + y^2)$$
$$W(p_{5,5}) = y(x+y)(-x^3 + 6x^2y + 4xy^2 + y^3)$$
$$W(p_{6,6}) = y(x+y)(16x^4 + 10x^3y + 10x^2y^2 + 5xy^3 + y^4)$$
In general, $W(p_{n,n}) = y(x+y)Q_n(x,y)$, where $Q_n$ is a polynomial in $x$ and
$y$ of degree $n-2$.  $W(p_{n,n})$ has a factor of $(x+y)^2$ if and only if $Q_n$
has a factor of $(x+y)$; i.e. if $Q_n(-y,y) = 0$. From the Lemma~\ref{L:circle}, we know:
$$Q_n = y(3x+y)(x+y)^{n-4} + 3x^2Q_{n-2} - 2x^3Q_{n-3}$$
In particular, for $n>6$ we have $Q_n(-y,y) = 3y^2Q_{n-2}(-y,y) + 2y^3Q_{n-3}(-y,y)$.
And from above, we see that $Q_4(-y,y) = 4y^2$, $Q_5(-y,y) = 4y^3$,
$Q_6(-y,y) = 12y^4$; in particular, the coefficients are all {\it positive}.  By
induction, the coefficient $c_n$ of $Q_n(-y,y) = c_ny^{n-2}$ is monotonically
increasing and always positive, and therefore can never be 0.  We conclude that
$Q_n$ does {\it not} have a factor of $(x+y)$, which completes the lemma. $\Box$

So now we can prove our final proposition:
\begin{prop} \label{P:independent}
$\{p_n, p_{n,k}\}_{k=4}^n$ is an independent set for all $n$.
\end{prop}
{\sc Proof:}  We will prove this by showing that $\{W(p_n), W(p_{n,k})\}_{k=4}^n$
is an independent set.  Again, we use induction.  The base case is clear from
our computations.  Assume that $\{W(p_{n-1}), W(p_{n-1,k})\}_{k=4}^{n-1}$ is an
independent set.  Then Lemma~\ref{L:weights} implies that $\{W(p_n), W(p_{n,k})\}_{k=4}^{n-1}$
is an independent set.  It remains only to show that $W(p_{n,n})$ is independent
from the other elements.  But we have just seen that the other elements are all
divisible by $(x+y)^2$, whereas $W(p_{n,n})$ is not (by Lemma~\ref{L:factor}), so it cannot
be a linear combination of the others.  Hence $\{W(p_n), W(p_{n,k})\}_{k=4}^n$
is independent for all $n$, so $\{p_n, p_{n,k}\}_{k=4}^n$ must also be independent
for all $n$.  $\Box$

This completes the proof of Theorem~\ref{T:loop_algebra}.

\section{Questions and Conjectures} \label{questions}
Although we know from studying mutant knots that the Intersection Graph
Conjecture fails in general (see \cite{cd}, \cite{le} and \cite{mc}), there are still many questions
left to be asked.  It is still unknown how badly the conjecture fails.  In order to determine
exactly how useful intersection graphs are in the study of chord diagrams and
finite type invariants, we would need to answer the following question:
\begin{quest} \label{Q:kernel}
What is the kernel of the map $\Gamma$, in each degree, from the space of
chord diagrams modulo the 4-term relation to the space of intersection graphs,
modulo the images of all 4-term relations?
\end{quest}
The results of Chmutov et. al.~\cite{cdl2} and Theorem~\ref{T:IGCloops} show that this kernel
is trivial if the map is restricted to the space of tree and loop diagrams.
Chmutov et. al.~\cite{cdl1} have shown via computer calculations that the Intersection Graph
Conjecture holds for chord diagrams of degree $\leq 8$.
\begin{prop} \label{P:trivial}
The kernel of $\Gamma$ is trivial when restricted to chord diagrams of degree $\leq 8$.
\end{prop}
{\sc Proof:}  Assume that a linear combination $\sum{k_iD_i}$ of chord diagrams of degree $n \leq 8$
is in the kernel of $\Gamma$.  Then the image $\sum{k_i\Gamma(D_i)}$ is trivial modulo the images of
all 4-term relations.  So by adding the images of some 4-term relations to the linear combination of
intersection graphs, we can cancel all of the graphs.  Since the Intersection Graph Conjecture holds
for diagrams of degree $\leq 8$, each graph that we add corresponds to a unique chord diagram; so
adding the images of the 4-term relations to the combination of intersection graphs corresponds to
adding a unique set of 4-term relations to the linear combination of chord diagrams.  Since all the
graphs in the resulting combination of intersection graphs cancel, so must the corresponding chord
diagrams in the combination of chord diagrams (again, because each graph corresponds to a unique chord
diagram).  Therefore, $\sum{k_iD_i}$ is trivial modulo the 4-term relation, and we are done.  $\Box$

However, by Morton and Cromwell~\cite{mc}, the kernel is known to be non-
trivial in degree 11.  Nothing else is known; in particular, we would like
to know if the kernel is trivial in degrees 9 and 10.

In addition, while the kernel is known to be non-trivial in degree 11, no-one
has exhibited two explicit inequivalent chord diagrams of degree 11 which have
the same intersection graph.  However, we can glean some obvious possibilities
from Morton and Cromwell~\cite{mc}.  Since they show that the Conway and Kinoshita-Terasaka knots
can be distinguished by a finite-type invariant of type 11, we can begin by
looking at the planar projections of these knots, which are then singular
knots with 11 double points.  The chord diagrams for these singular knots are:
$$\psfig{file=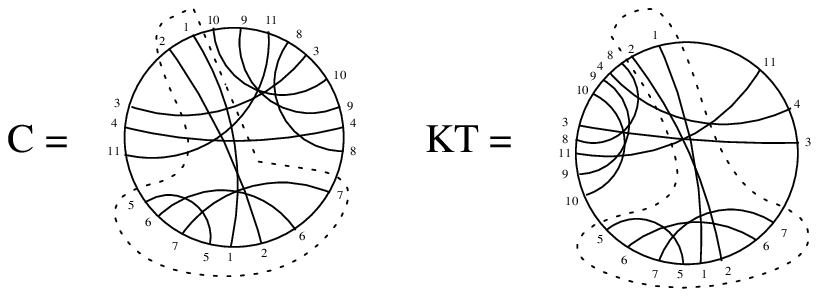}$$
These chord diagrams differ just by rotating one share about another (the fixed
share is the one enclosed by the dashed line), so they
have the same intersection graph.  The labeling of the chords in the chord
diagrams corresponds to the labeling of the chords in the intersection graph:
$$\psfig{file=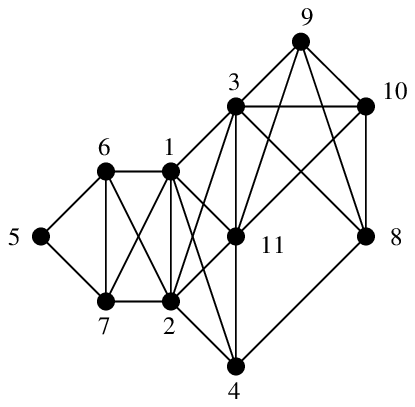}$$
We can then conjecture that:
\begin{conj} \label{C:C=KT}
The chord diagrams C and KT are not equivalent.
\end{conj}

We can also ask how far the approach of this paper, looking at the rank of
the fundamental group of the intersection graph, rather than at its degree,
can be extended.  The fact that the conjecture fails at degree 11 means that
there is a counterexample with at most $\left(\matrix{11 \cr 3}\right) = 165$
loops.  The counterexample proposed above has 15 loops.  However, it does not
seem that the arguments used in this paper can be easily extended to the case
of diagrams with two loops, so perhaps there is counterexample there.
\begin{quest} \label{Q:2loop}
Is there a counterexample to the IGC involving 2-loop diagrams?
\end{quest}

A first step towards answering Question~\ref{Q:kernel} would be to determine in general
the group of elementary transformations which can be performed on a chord
diagram without changing its intersection graph (i.e., find some set of
generators for this group).  It seems likely that we would have to describe
such a set of generators in terms of shares:
\begin{conj} \label{C:elementary}
The group of elementary transformations is generated by transformations of the
following two types:
\begin{itemize}
    \item Reflecting a share across another share:
$$\psfig{file=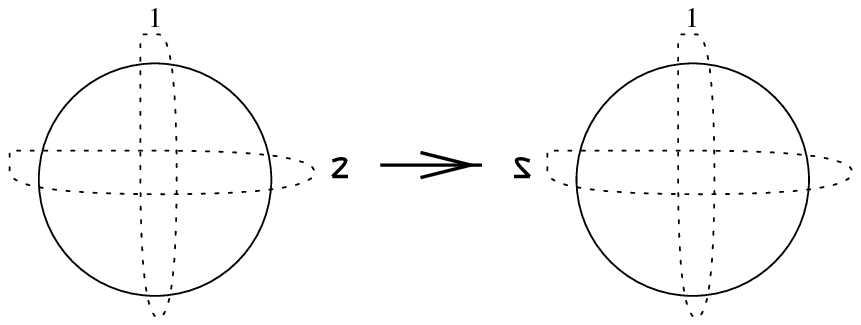}$$
    \item Turning a share upside down:$$\psfig{file=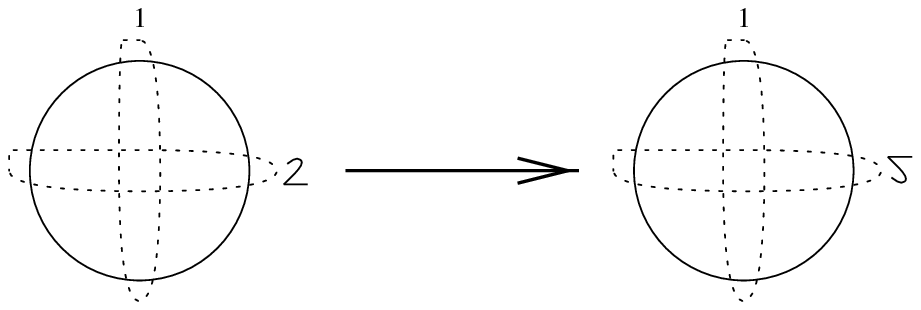}$$
\end{itemize}
\end{conj}
These transformations will certainly generate all transpositions, and hence
all permutations, of shares along a chord (or another share); so they generate
all the elementary transformations used in \cite{cdl2} and in this paper.

A final question concerns the size of the equivalence classes of intersection
graphs under the relations induced by $\Gamma$.  The results of Chmutov et. al.~\cite{cdl2}
and this paper have demonstrated that there are different chord diagrams
(individual chord diagrams, not linear combinations) which are equivalent
modulo the 4-term relation - namely, tree and loop diagrams sharing the same
intersection graph.  But these results do not help us answer the analogous
question in the space of intersection graphs:
\begin{quest} \label{Q:equivalent_graphs}
Are there two different intersection graphs (individual graphs, not linear
combinations) which are equivalent modulo the 4-term relations (i.e. the
relations induced by the 4-term relations via $\Gamma$)?
\end{quest}
It seems likely that the answer is ``Yes,'' as it is for chord diagrams, but
it would be interesting to have an explicit example.  It would be even better
to discover the ``elementary transformations'' between equivalent intersection
graphs.

Intersection graphs may still have much to offer us as a tool for studying the
space of chord diagrams, but there is still a lot of work to be done before
we can exploit them fully.
\section{Acknowledgements}
I wish to thank the following people for many valuable conversations and criticisms:
Paul Melvin, Sergei Chmutov, David Gay, Andrew Lewis, Robert Schneiderman and especially
my advisor Robion Kirby.  I also wish to thank the anonymous reviewer who pointed out
the fatal flaw in my original proof of Theorem~\ref{T:loop_algebra}, and the second reviewer
for his or her many helpful comments.  Finally, I want to acknowledge the support of UC
Berkeley, where this work was done while I was a graduate student.
\end{document}